\documentclass[
aps,%
11pt,%
final,%
notitlepage,%
oneside,%
onecolumn,%
nobibnotes,%
nofootinbib,% 
superscriptaddress,%
noshowpacs,%
centertags]%
{revtex4}

\usepackage[letterpaper, nohead,dvips]{geometry}
\usepackage{amsmath}
\usepackage{amsfonts}
\usepackage{amssymb}
\usepackage{enumerate}
\usepackage{amsthm}
\usepackage{setspace}
\usepackage[all,cmtip]{xy}
\usepackage{graphicx}
\usepackage{epstopdf}
\usepackage{MnSymbol}

\newtheorem{theorem}{Theorem}
\newtheorem{lemma}{Lemma}

\newcommand{\curlyd}{\partial}

\newcommand {\Ttensor}{S}
\newcommand{\ctensor}{ {\mathcal C}}
\newcommand{\btensor}{{\mathcal B}}

\begin{document}

\title{On the Method of  Interconnection and Damping Assignment Passivity-Based Control for the Stabilization of Mechanical Systems}

\author{\firstname{Dong Eui} \surname{Chang}}
\email{dechang@uwaterloo.ca}
\affiliation{Department of Applied Mathematics, University of Waterloo,
200 University Avenue West, Waterloo, ON   N2L 3G1, Canada }

\begin{abstract}
Interconnection and damping assignment passivity-based control (IDA-PBC) is an excellent method to stabilize mechanical systems in the Hamiltonian formalism. In this paper, several improvements are made on
the IDA-PBC method. The skew-symmetric interconnection submatrix  in the conventional form of IDA-PBC is shown to have some redundancy for systems with the number of degrees of freedom greater than two,  containing unnecessary components that do not contribute to the dynamics. To completely remove this redundancy,
the use of quadratic gyroscopic forces is proposed in place of the skew-symmetric interconnection submatrix.  Reduction of the number of  matching partial differential equations in IDA-PBC and simplification of
the structure of the matching partial differential equations  are achieved by eliminating the gyroscopic force from the matching partial differential equations.   In addition,  easily verifiable criteria are provided for Lyapunov/exponential stabilizability by IDA-PBC for all linear controlled Hamiltonian  systems with  arbitrary degrees of underactuation and for all nonlinear controlled Hamiltonian  systems with one degree of underactuation.  A general design procedure for IDA-PBC is given and  illustrated with examples.  The duality of the new IDA-PBC method to the method of controlled Lagrangians is discussed. This paper renders the IDA-PBC method as powerful as the controlled Lagrangian method.

\vspace{3mm}\noindent

\vspace{3mm}\noindent

\vspace{3mm}\noindent 
 \underline{Published in {\it Regular and Chaotic Dynamics}, 19 (5), 556 -- 575, September  2014.}
\end{abstract}

\maketitle

\section{Introduction}

Mechanical systems are ubiquitous in nature and engineering, and there have been many studies on modeling motions of insects, animals, fish and humans in the framework of mechanics and control. For example, a model for clock-actuated legged locomotion  of human and biologically-inspired robots is proposed and analyzed in \cite{SpHo07}. A control-theoretic strategy for human walking gait assistance is suggested with a  biped model to lessen the perceived weight of a patient's center of mass through a  robotic angle-foot orthosis with one actuated degree-of-freedom \cite{GrBrSp10}. A time-scaling control law is developed and applied to two passive-dynamic bipeds: a compass-gait biped and a simple biped with torso  \cite{HoLeSp07}.  One of the main objectives in all these studies is to study stability and stabilization of the motion of the system.

The energy shaping method stands out among the methods for the stabilization of mechanical systems since it preserves the mechanical structure, provides a systematic procedure for constructing control laws, and yields a large region of stability. The idea of this method is simple. Given an unstable mechanical system, one  transforms it via feedback to a stable mechanical system whose total energy function obtains a minimum value at the equilibrium of interest and then injects a dissipative feedback force to obtain asymptotic stability of the equilibrium point. In this process of transformation, the original total energy function with a saddle-type critical point at the equilibrium point is transformed to a new total energy function with a minimum value at the equilibrium point.  In order to find such a  new mechanical system with a stable energy function, one has to solve partial differential equations (PDEs) for the mass matrix and the potential function of the new mechanical system. The PDEs for the new mass matrix are called kinetic matching conditions or kinetic matching PDEs, and the PDEs for the new potential function are called potential matching conditions or potential matching PDEs. Hence, understanding the structure of the matching PDEs and their solvability is crucial in the application of the energy shaping method.

To help  readers grasp the idea quickly, let us give some examples of energy shaping.  Consider the following one-degree-of-freedom mechanical system with  control $u$:
\begin{align*}
\dot q^1 = p_1; \quad
\dot p_1 = q^1 + u,
\end{align*}
where $q^1, p_1, u \in \mathbb R$. It can be written in Hamiltonian form with the control force $u$ as follows:
\[
\begin{bmatrix}
\dot q^1 \\ \dot p_1
\end{bmatrix} = \begin{bmatrix}
0 & 1 \\ -1 & 0
\end{bmatrix}
\begin{bmatrix}
\frac{\partial H}{\partial q^1} \\ 
\frac{\partial H}{\partial p_1}
\end{bmatrix} +\begin{bmatrix}
0 \\ 1
\end{bmatrix}u
\]
where 
\begin{equation}\label{intro:H}
H(q^1,p_1) = \frac{1}{2}(p_1)^2 - \frac{1}{2}(q^1)^2
\end{equation}
is the Hamiltonian (or total energy) function. The equilibrium point $(q^1,p_1) = (0,0)$ is unstable in the uncontrolled dynamics with $u=0$ and it is a saddle point of $H$. Apply the feedback control 
\begin{equation}\label{intro:u}
u = -2q^1 -b p_1
\end{equation}
with $b>0$ 
 to transform the system to the following closed-loop system
\[
\dot q^1 = p_1; \quad \dot p_1 = -q^1 -b p_1
\]
which can be put in the following form
\begin{equation}\label{first:mat}
\begin{bmatrix}
\dot q^1 \\ \dot p_1
\end{bmatrix} = \begin{bmatrix}
0 & 1 \\ -1 &-b
\end{bmatrix}
\begin{bmatrix}
\frac{\partial \widehat H}{\partial q^1} \\ 
\frac{\partial \widehat H}{\partial p_1}
\end{bmatrix} 
\end{equation}
where 
\begin{equation}\label{intro:Hhat}
\widehat H (q^1,p_1) =  \frac{1}{2}(p_1)^2 + \frac{1}{2}(q^1)^2.
\end{equation}
The new Hamiltonian $\widehat H$ obtains its minimum value at the equilibrium point $(q^1,p_1) = (0,0)$ and its time-derivative along the trajectory of the closed-loop system is computed as
\begin{equation}\label{intro:H:dec}
\frac{d \widehat H}{dt} = -b (p_1)^2 \leq 0,
\end{equation}
which implies that $(q^1,p_1) = (0,0)$ is Lyapunov stable.  By computing the eigenvalues of the closed-loop dynamics, one can show that $(q^1,p_1) = (0,0)$ is asymptotically stable. Notice that the first term $-2q^1$ of the control $u$ in (\ref{intro:u}) shapes the original Hamiltonian $H$ to the new one $\widehat H$ by altering the potential function from $-\frac{1}{2}(q^1)^2$ to $\frac{1}{2}(q^1)^2$ so that the equilibrium point $(0,0)$ becomes the minimum point of the new Hamiltonian function.  The second term $-b p_1$ of the control $u$ in (\ref{intro:u}) introduces dissipation to the closed-loop system such that the Hamiltonian $\widehat H$ becomes non-increasing in the closed-loop system as shown in (\ref{intro:H:dec}) to achieve asymptotic stability. 

Sometimes, potential energy shaping alone is not enough, so it is necessary to do kinetic energy shaping as well. Consider the following system:
\begin{align*}
\dot q^1 &= p_1,\\
\dot q^2 &= p_2, \\
\dot p_1 &= q^1 - q^2,\\
\dot p_2 & = -q^1 + u,
\end{align*}
where $(q^1,q^2) \in \mathbb R^2$ is the position vector; $(p_1, p_2) \in \mathbb R^2$ is the momentum vector; and $u \in \mathbb R$ is the control. It is easy to verify that the equilibrium point $(q^1,q^2,p_1,p_2) = (0,0,0,0)$ is unstable when $u=0$.  The above system can be put in controlled Hamiltonian form as follows:
\begin{equation}\label{int:old:eq}
\begin{bmatrix}
\dot q^1 \\ \dot q^2 \\ \dot p_1 \\ \dot p_2
\end{bmatrix} = \begin{bmatrix}
0 & 0 & 1 & 0 \\
0 & 0 & 0 & 1 \\
-1 & 0 & 0 & 0\\
0 & -1 & 0 & 0
\end{bmatrix}\begin{bmatrix}
\frac{\partial H}{\partial q^1} \\ \frac{\partial H}{\partial q^2} \\ \frac{\partial H}{\partial p_1} \\ \frac{\partial H}{\partial p_2}
\end{bmatrix} + \begin{bmatrix}
0 \\ 0 \\ 0 \\1
\end{bmatrix}u,
\end{equation}
 where
\begin{equation}\label{intro:H0}
H = \frac{1}{2}(p_1)^2 + \frac{1}{2}(p_2)^2 - \frac{1}{2}(q^1)^2 + q^1 q^2
\end{equation}
is the Hamiltonian. Notice that $H$ does not obtain a minimum value at the origin, being consistent with the instability of the origin.  Apply the following feedback control:
\begin{equation}\label{intro:u:2}
u = 7 q^1 - 4 q^2 + 5p_1 - 2p_2
\end{equation}
which transforms the given system to the following closed-loop system
\begin{equation}\label{int:new:eq}
\begin{bmatrix}
\dot q^1 \\ \dot q^2 \\ \dot p_1 \\ \dot p_2
\end{bmatrix} = \begin{bmatrix}
0 & 0 & 2 & 5 \\
0 & 0 & 5 & 13 \\
-2 & -5 & 0 & 0\\
-5 & -13 & 0 & -1
\end{bmatrix}\begin{bmatrix}
\frac{\partial \widehat H}{\partial q^1} \\ \frac{\partial \widehat H}{\partial q^2} \\ \frac{\partial \widehat H}{\partial p_1} \\ \frac{\partial \widehat H}{\partial p_2}
\end{bmatrix} 
\end{equation}
 where
\begin{equation}\label{intro:H1}
\widehat H = \frac{13}{2}(p_1)^2 -5p_1p_2 + (p_2)^2 + \frac{17}{2}(q^1)^2 -7q^1q^2 + \frac{3}{2}(q^2)^2.
\end{equation}
The new Hamiltonian $\widehat H$ has its minimum value at the origin and its time derivative along the trajectory of the closed-loop system satisfies
\[
\frac{d \widehat H}{dt} = -(5p_1-2p_2)^2 \leq 0,
\]
proving Lyapunov stability of the origin. One can further show that the origin is asymptotically stable. Notice that the term $7 q^1 - 4 q^2$ in the control $u$ in  (\ref{intro:u:2}) shapes the unstable Hamiltonian function $H$ to the stable  Hamiltonian function $\widehat H$ changing both the kinetic energy and the potential energy. The term $5p_1 - 2p_2$ in the control $u$ in  (\ref{intro:u:2}) injects damping for asymptotic stability.  Notice also that the control $u$  modifies the Hamiltonian structure from the canonical structure in (\ref{int:old:eq}) to the non-canonical structure in (\ref{int:new:eq}), where the matrix in (\ref{int:new:eq}) is the sum of a symplectic matrix and a negative semidefinite symmetric damping matrix. 

Energy shaping for stabilization is not always possible. Consider the following system:
\begin{align*}
\dot q^1 &= p_1,\\
\dot q^2 &= p_2, \\
\dot p_1 &= -q^1,\\
\dot p_2 & = q^2 + u,
\end{align*}
where $(q^1,q^2) \in \mathbb R^2$ is the position vector; $(p_1, p_2) \in \mathbb R^2$ is the momentum vector; and $u \in \mathbb R$ is the control.  Notice that the sub-system
\begin{align*}
\dot q^1 = p_1; \quad
\dot p_1 = -q^1
\end{align*}
is exponentially unstable and decoupled from the rest of the system. Hence, it is impossible to apply a feedback control to shape the Hamiltonian function such that the new Hamiltonian function obtains its minimum value at the origin and its time derivative is non-positive along the trajectory of the closed-loop system. 

In  the first two examples  above, energy shaping controllers were given without derivation, but  it is important to have  a systematic procedure for energy shaping. The third example shows the necessity of an easily verifiable   criterion for stabilizability by the energy shaping method. These issues are addressed later in this paper.

Let us briefly review the history of development of the energy shaping method. The idea of  potential energy shaping for the stabilization of  mechanical systems dates back to \cite{ArMi83}. The notion of kinetic energy shaping first appeared in \cite{BlKrMaAl92}.  The idea of total energy shaping was introduced in \cite{BlLeMa97,BlLeMa00}, and has then been actively developed in \cite{BlChLeMa01,BlLeMa01:IJRNC,Ch08:IFAC,Ch10:SIAM,Ch10:TAC,Ch12:IJRNC,ChBlLeMaWo02,GaLeMa08,NgChLa13,Ze02}.
There are two approaches to energy shaping: the Lagrangian approach and the Hamiltonian approach. The energy shaping method is called the method of controlled Lagrangians on the Lagrangian side and the method of interconnection and damping assignment passivity-based control (IDA-PBC) on the Hamiltonian side. It has been proven that the two approaches are equivalent \cite{Ch08:IFAC}.  Chang \cite{Ch10:SIAM} has then improved the Lagrangian method  by completely characterizing quadratic gyroscopic forces, reducing the number of matching conditions  for energy shaping, and finding necessary and sufficient conditions for stabilizability by energy shaping for the class of all linear mechanical systems with arbitrary degrees of underactuation and the class of all nonlinear mechanical systems with underactuation degree one. In contrast to these developments on the Lagrangian side, there has been a lag in  development on the Hamiltonian side except for \cite{Ch10:Med}.

In this paper, we make improvement of the method of IDA-PBC in several ways so that the IDA-PBC method becomes as powerful as its Lagrangian counterpart. First, we propose a new form of IDA-PBC by introducing a quadratic gyroscopic force.  The skew-symmetric interconnection matrix in the conventional form of IDA-PBC \cite{AcOrAsMa05,OrGa04,RoDoOr13} is shown to have unnecessary components that do not appear in the system dynamics with the number of degrees of freedom greater than two. To remove the unnecessary components, the use of quadratic gyroscopic force is proposed in place of the skew-symmetric interconnection matrix. Second, the kinetic matching conditions get decomposed into two sets: one containing no components of the gyroscopic force and the other containing components of the gyroscopic force. The first set constitutes a new reduced set of kinetic matching PDEs and the second is used to   algebraically define  the gyroscopic force.  This decomposition substantially reduces the number of kinetic matching PDEs. For example, if the degree of underactuation is one,  there is only one  kinetic matching PDE, irrespective of the number of degrees of freedom.  Moreover,  when the co-distribution spanned by actuation (or control) co-vector fields is integrable,   the new reduced matching PDEs contain a smaller number of entries of the unknown mass matrix, which was never discovered in the past for IDA-PBC.  Third, we provide  necessary and sufficient conditions for stabilizability by IDA-PBC for the class of all linear mechanical systems with an arbitrary degree of  underactuation and for the class of all nonlinear mechanical systems with one degree of underactuation.   These conditions can be easily verified in advance before solving the matching PDEs.
  Fourth, a step-by-step synthesis procedure with IDA-PBC is provided and the main results are illustrated with examples. Lastly, the duality of the new method of  IDA-PBC  to the method of controlled Lagrangians is discussed.

\section{Main Results}

We  review some basic notions on tensors, derive important lemmas, present a new form of IDA-PBC and  compare it with the conventional form.   Although both forms are shown to be theoretically equivalent for use in IDA-PBC,  the new one has the advantage that it contains fewer components.  Using the new form, we derive matching conditions for energy shaping and  decompose them into two parts:  one without any gyroscopic terms and the other with gyroscopic terms.  This decomposition  allows  us to reduce the number of matching PDEs. We  discuss conditions for stabilizability by IDA-PBC that can be verified easily, and illustrate the results with  examples. The duality of the IDA-PBC method to the controlled Lagrangian method is discussed.

\subsection{Notation}

Let $V$ be an $n$-dimensional real vector space and $V^*$ its dual space.
Let $\{e_1, \ldots, e_n\}$ be a basis of $V$, and $\{e^1, \ldots, e^n\}$ its dual basis such that $\langle e^i, e_j\rangle = \delta^i_j$, where $\langle, \rangle$ is the canonical pairing between dual vectors and vectors, and $\delta^i_j$ is the Kronecker delta. The tensor product of two vector spaces $V$ and $W$ is denoted by $V \otimes W$, and each element of $V \otimes W$ is a linear combination of elements of the form $v \otimes w$, where $v \in V$ and $w\in W$. The $r$-fold tensor product of a vector space $V$ is denoted by $V^{\otimes r}$ or $\otimes^rV$.
An $(r,s)$-tensor  $T$  on $V$ is written as
\[
T = T^{i_1 \cdots i_r}_{j_1 \cdots j_s} e_{i_1} \otimes \cdots \otimes e_{i_r} \otimes e^{j_1} \otimes \cdots \otimes e^{j_s},
\]
where the Einstein summation convention is enforced.  The contraction of  an $(r,0)$-tensor  $S = S^{i_1 \cdots i_r} e_{i_1} \otimes \cdots \otimes e_{i_r}$ and a $(0,s)$-tensor $T = T_{j_1 \cdots j_s} e^{j_1} \otimes \cdots \otimes e^{j_s}$ with $r \leq s$,  is  defined and denoted by
\[
S\righthalfcup T =  S^{i_1 \cdots i_r} T_{i_1 \cdots i_r i_{r+1} \cdots i_s} e^{i_{r+1}} \otimes \cdots \otimes e^{i_s}.
\]
 One can identify each $(0,s)$-tensor $T$ with a multi-linear map $\tilde T : V \times \cdots \times V \rightarrow \mathbb R$ defined by
\begin{align*}
\tilde T(v_1, \cdots, v_s) &= v_1 \otimes \cdots \otimes v_s \righthalfcup T\\
&= v_s \righthalfcup \cdots \righthalfcup v_1 \righthalfcup T.
\end{align*}
Define the symmetrization operator $\operatorname{Sym} : V^{\otimes r} \rightarrow V^{\otimes r}$ by
\[
\operatorname{Sym} (v_1 \otimes \cdots \otimes v_r) = \frac{1}{r!}\sum_{\sigma \in {\rm S}_r} v_{\sigma (1)}\otimes \cdots \otimes v_{\sigma (r)},
\]
where  ${\rm S}_r$ is the symmetric group of $\{1, \ldots, r\}$ and it is understood that $\operatorname{Sym}$ is linearly extended to   $V^{\otimes r}$.

A (0, 3)-tensor $C = C_{ijk}e^i \otimes e^j \otimes e^k$ on $V$ is said to be gyroscopic if it satisfies
\begin{align}
&C(u,v,w)  = C(v, u,w) \quad \forall u,v,w \in V, \label{C:symmetry:intrinsic}\\
& C(u,v,w) + C(v,w,u) + C(w,u,v) = 0 \quad \forall u,v,w \in V.\label{C:Jacobi:intrinsic}
\end{align}
 In coordinates, (\ref{C:symmetry:intrinsic}) and (\ref{C:Jacobi:intrinsic}) are equivalent, respectively, to the following
\begin{align}
&C_{ijk} = C_{jik}, \label{C:symmetry}\\
&C_{ijk} + C_{jki} + C_{kij} = 0\label{C:Jacobi}
\end{align}
for all $i, j, k$.
Define
\begin{align*}
\ctensor (V) &= \{ C \in V^{*\otimes 3} \mid \textup{properties (\ref{C:symmetry:intrinsic}) and (\ref{C:Jacobi:intrinsic}) hold} \} \\
&= \{ C_{ijk}e^i \otimes e^j \otimes e^k \mid\textup{properties (\ref{C:symmetry}) and (\ref{C:Jacobi}) hold}\}.
\end{align*}
Define
\begin{align*}
\btensor(V)& = \{ B \in V^{*\otimes 3} \mid B(u,v,w) = -B(u,w,v) \forall u,v,w \in V\}\\
&= \{B_{ijk}e^i \otimes e^j \otimes e^k \mid B_{ijk} = -B_{ikj} \}.
\end{align*}
For a manifold $M$, $\btensor (M)$ and $\ctensor (M)$ denote the sets of the tensor fields on $M$ such that at each point $q\in M$, $\btensor (M)_q$ and $\ctensor (M)_q$ are equal to $\btensor (T_qM)$ and $\ctensor (T_qM)$, respectively, where $T_qM$ denotes the tangent space to $M$ at $q$.

Below are four lemmas, all of which can be skipped  until they are referred to later in this paper. 

\begin{lemma}\label{dim:BC}
Let $n = \dim V$. Then, $\dim \btensor (V) = \frac{n^2(n-1)}{2}$ and $\dim \ctensor (V) = \frac{n(n^2-1)}{3}$.
\begin{proof}
Due to the skew-symmetry property $B_{ijk} = -B_{ikj}$, we have $\dim \btensor (V) = \frac{n^2(n-1)}{2}$. To compute $\dim \ctensor (V)$, we consider three cases: 1) all three identical indices; 2) two identical indices and one distinct index; and 3) three distinct indices. For case 1, we only have $C_{iii} = 0$ by (\ref{C:Jacobi}). For case 2, we have $C_{iki} = C_{kii} = -\frac{1}{2}C_{iik}$  for all $i \neq k$ by (\ref{C:symmetry}) and (\ref{C:Jacobi}). Hence, the number of independent components in this case is $n(n-1)$. For case 3, the number of independent components is $\frac{n(n-1)(n-2)}{3}$ by (\ref{C:symmetry}) and (\ref{C:Jacobi}). Hence, $\dim \ctensor (V) = 0+n(n-1) + \frac{n(n-1)(n-2)}{3} = \frac{n(n^2-1)}{3}$.
\end{proof}
\end{lemma}

\begin{lemma} \label{lemma:BC}
Let $\psi : \btensor (V) \rightarrow V^{*\otimes 3}$ be a linear map defined by
\begin{equation}\label{def:psi}
\psi(B_{ijk}e^i \otimes e^j \otimes e^k) = \frac{1}{2} (B_{ijk} + B_{jik})e^i \otimes e^j \otimes e^k
\end{equation}
for  $B_{ijk}e^i \otimes e^j \otimes e^k \in \btensor (V)$.
Then, $\psi (\btensor (V)) = \ctensor (V)$, i.e., the image of $\btensor (V)$ under the map $\psi$ equals $\ctensor (V)$. Moreover, $\dim \operatorname{ker}\psi = \frac{n(n-1)(n-2)}{6}$, where $n=\dim V$.
\begin{proof}
First, we show that $\psi (\btensor (V)) \subset \ctensor (V)$.
Choose any $B = B_{ijk}e^i \otimes e^j \otimes e^k \in B(V)$. Let $C = \psi (B)$. In coordinates, $C_{ijk} =  \frac{1}{2} (B_{ijk} + B_{jik})$, where $C = C_{ijk} e^i \otimes e^j \otimes e^k$.  Then, $C_{ijk} = \frac{1}{2} (B_{ijk} + B_{jik}) = \frac{1}{2} (B_{jik} + B_{ijk}) = C_{jik}$, satisfying (\ref{C:symmetry}). Also, $C_{ijk} + C_{jki} + C_{kij} = \frac{1}{2} (B_{ijk} + B_{jik}) + \frac{1}{2}  (B_{jki} + B_{kji}) + \frac{1}{2} (B_{kij} + B_{ikj} ) = \frac{1}{2}  (B_{ijk} + B_{jik} - B_{jik} + B_{kji} - B_{kji} - B_{ijk})=0$, satisfying (\ref{C:Jacobi}).
Hence, $\psi(\btensor (V)) \subset \ctensor (V)$.

We now show $\ctensor (V) \subset \psi(\btensor (V))$. Choose any $C = C_{ijk} e^i \otimes e^j \otimes e^k \in \ctensor (V)$ such that $C_{ijk}$ satisfies (\ref{C:symmetry}) and (\ref{C:Jacobi}). It  suffices to choose a tensor $B = B_{ijk}e^i \otimes e^j \otimes e^k \in \btensor (V)$ such that $\psi(B) = C$. Let $B_{iii} = 0$ for all $i$. For all $i\neq k$, let $B_{iik} = C_{iik}$, $B_{iki} = -C_{iik}$ and $B_{kii} = 0$. For all $i<j<k$, let
\begin{align*}
B_{ijk} &= 2C_{ijk}, \quad B_{kij} = -2C_{jki}, \quad  B_{jki} = 0,\\
B_{ikj} &= -2C_{ijk}, \quad B_{kji} = 2C_{jki}, \quad B_{jik} = 0.
\end{align*}
Then, it is easy to show  $B \in \btensor (V)$  and $\psi(B) = C$. Hence, $\ctensor (V) \subset \psi(\btensor (V))$. Since $\psi (\btensor (V)) \subset \ctensor (V)$ and $\ctensor (V) \subset \psi(\btensor (V))$, it follows that $\psi (\btensor (V)) = \ctensor (V)$.

Since $\psi$ is onto,  $\dim \operatorname{ker}\psi = \dim \btensor (V) - \dim \ctensor(V) =  \frac{n(n-1)(n-2)}{6}$ by Lemma \ref{dim:BC}.
\end{proof}
\end{lemma}

\begin{lemma}\label{first:lemma}
1. For a $(0,s)$-tensor $T$ on a vector space $V$, $T( v, \ldots, v) = 0$ for all $v \in V$ if and only if $\operatorname{Sym}(T) = 0$.

2. Suppose that a $(0,3)$-tensor $C$ on $V$ satisfies (\ref{C:symmetry:intrinsic}). Then, $\operatorname{Sym}(C) = 0$ if and only if
$C$ satisfies (\ref{C:Jacobi:intrinsic}).
\begin{proof}
(a) Trivial by definition of $\textrm{Sym}$.

(b) Let $C = C_{ijk} e^i \otimes e^j \otimes e^k$. Then, $0=\operatorname{Sym}(C)$ $\Leftrightarrow$ $0 = \frac{1}{6}( C_{ijk} + C_{jki} + C_{kij} + C_{jik} + C_{kji} + C_{ikj} ) = \frac{1}{3} (C_{ijk} + C_{jki} + C_{kij} ) $ $\Leftrightarrow$ $C$ satisfies (\ref{C:Jacobi:intrinsic}).
\end{proof}
\end{lemma}

\begin{lemma} \label{lemma:partial:gyro}
Let $ V = \tilde V \oplus \hat V$, and $\Ttensor = \Ttensor_{ijk} e^i \otimes e^j \otimes e^k$ be a $(0,3)$-tensor on $V$ that satisfies
 \begin{equation}\label{T:sym}
 \Ttensor (u,v,w) =\Ttensor (v,u,w) \quad \forall u,v,w\in  V.
 \end{equation}
Then, the following are equivalent:
\begin{description}
\item [(a)] There exists a $(0,3)$-tensor  $C = C_{ijk}e^i \otimes e^j \otimes e^k$ that satisfies (\ref{C:symmetry:intrinsic}), (\ref{C:Jacobi:intrinsic}) and
\begin{equation}\label{CT}
C(v_1,v_2,\tilde v) = \Ttensor (v_1,v_2,\tilde v) \,\,\forall v_1,v_2\in V, \tilde v\in \tilde V.
\end{equation}

\item[(b)] $\Ttensor (\tilde v_1,\tilde v_2,\tilde v_3) + \Ttensor (\tilde v_2,\tilde v_3,\tilde v_1) + \Ttensor (\tilde v_3,\tilde v_1,\tilde v_2) =0 \quad \forall \tilde v_1, \tilde v_2, \tilde v_3\in \tilde V$.
\end{description}
\begin{proof}
(a) $\Rightarrow $ (b): Trivial.

(b) $\Rightarrow $ (a): We prove this in coordinates since it will be conveniently used later. Let $\dim V = n$, $\dim \tilde V = n -m$ and $\dim \hat V = m$. Use the following three groups of indices: $i,j,k =1, \ldots, n$; $\alpha, \beta, \gamma = 1, \ldots, n -m $; and $a, b, c =(n-m+1), \ldots, n$. Suppose $V = \operatorname{span}\{ e_i\}$, $\tilde V = \operatorname{span} \{e_\alpha\}$ and $\hat V = \operatorname{span}\{e_a\}$. Define $C_{ijk}$ as follows:
\begin{align}
C_{ij\alpha} &= \Ttensor_{ij\alpha},\label{C1}\\
C_{\alpha \beta a} &= -\Ttensor_{\beta a \alpha} - \Ttensor_{a \alpha \beta},\label{C2}\\
C_{\alpha ab} &= C_{b \alpha a  } = -\frac{1}{2}\Ttensor_{ab\alpha},\label{C3}\\
C_{abc} &=0.\label{C4}
\end{align}
It is then easy to show that $C_{ijk}$ satisfies (\ref{C:symmetry}), (\ref{C:Jacobi}) and (\ref{CT}). Instead of $C_{abc} = 0$, one can alternatively choose arbitrary $C_{abc}$ such that $C_{abc} = C_{bac}$ and $C_{abc} + C_{bca} + C_{cab} = 0$.
\end{proof}
\end{lemma}

\subsection{New Form of IDA-PBC for Mechanical Systems}
\label{section:new:form}

Let $\{e_1, \ldots, e_n\}$ be the standard basis of $\mathbb R^n$ and $\{e^1, \ldots, e^n\}$ its dual basis. For convenience, we identify $\mathbb R^n$ with its dual space $(\mathbb R^n)^*$, but we follow the convention that position vectors $q = q^ie_i$ are in $\mathbb R^n$ and momentum vectors $p = p_i e^i$ are in $(\mathbb R^n)^*$.  For the sake of simplicity we assume that every function is smooth.

Let us consider a controlled Hamiltonian system of the form
\begin{equation}\label{can:PBC}
\begin{bmatrix}
\dot q \\ \dot p
\end{bmatrix} = \begin{bmatrix}
0 & I_n \\
-I_n & 0
\end{bmatrix}
\begin{bmatrix}
\curlyd_q H \\ \curlyd_p H
\end{bmatrix} +
\begin{bmatrix}
0 \\ G(q)
\end{bmatrix}u
\end{equation}
where $q,p \in \mathbb R^n$, $u \in \mathbb R^m$, $m \leq n$, $I_n$ is the $n\times n$ identity matrix, $G(q)$ is an $n\times m$ matrix with $\operatorname{rank}G(q) = m$ for each $q$, and
\[
H(q,p) = \frac{1}{2}p^T M^{-1}(q)p + V(q)
\]
is the Hamiltonian function of the system, where $M(q)$ is an $n\times n$ positive definite symmetric mass matrix and $V(q)$ is a function called the potential function of the system. We say that this system has $n$ degrees of freedom and $(n-m)$ degrees of underactuation. Throughout this paper, we assume that $(q,p) = (0,0) \in \mathbb R^n \times \mathbb R^n$ is the equilibrium point of the system to be stabilized.
The controlled Hamiltonian system  (\ref{can:PBC}) shall be called linear if the matrices $M(q)$ and $G(q)$ are constant and the function $V(q)$ is a quadratic function of $q$. Otherwise, it shall be called nonlinear.

To achieve stabilization of the equilibrium point at the origin,  we have the objective of finding a feedback law $u = u(q,p)$ that transforms the system (\ref{can:PBC}) to the following desired form:
\begin{align}\label{desired:PBC:new}
\begin{bmatrix}
\dot q \\ \dot p
\end{bmatrix}  = \begin{bmatrix}
 0 & M^{-1}(q)\widehat M(q)\\
 -\widehat M(q)M^{-1}(q) & -G(q) K_{\rm d}(q) G^T  (q)
 \end{bmatrix} \begin{bmatrix}
\curlyd_q \widehat  H \\ \curlyd_p \widehat H
\end{bmatrix} + \begin{bmatrix}
0 \\
\widehat M^{-1} p \righthalfcup \widehat M^{-1} p \righthalfcup  C
\end{bmatrix},
\end{align}
where
\begin{equation}\label{hat:H}
\widehat H(q,p) = \frac{1}{2} p^T \widehat M^{-1}(q) p + \widehat V(q)
\end{equation}
is a desired Hamiltonian with an $n\times n$ positive definite symmetric mass matrix $\widehat M (q) $ and a function $\widehat V(q)$ having a non-degenerate minimum at $q=0$; $K_{\rm d} (q)$ is an $m\times m$ positive definite symmetric matrix, and $C = C_{ijk}(q)e^i \otimes e^j \otimes e^k$ is a $(0,3)$-tensor field that  pointwise satisfies (\ref{C:symmetry:intrinsic}) and (\ref{C:Jacobi:intrinsic}), i.e., $C \in \ctensor (\mathbb R^n)$. In (\ref{desired:PBC:new}), it is understood that
\begin{equation}\label{add:C}
\widehat M^{-1} p \righthalfcup \widehat M^{-1} p \righthalfcup  C = C_{ijk}\widehat M^{i\ell} \widehat M^{jr}p_\ell p_r e^k,
\end{equation}
where $\widehat M^{ij}$ denotes the $(i,j)$th entry of the inverse matrix $\widehat M^{-1}$ of $\widehat M$.
Without loss of generality, we have symmetrized $C_{ijk}$ with respect to its first two indices  (property (\ref{C:symmetry:intrinsic}) or (\ref{C:symmetry})) because $\widehat M^{-1} p$ appears quadratically in (\ref{add:C}). Notice  that $\operatorname{Sym}(C) = 0$ by statement 2 of Lemma \ref{first:lemma}. The function $\widehat V$ is called the potential energy of the system (\ref{desired:PBC:new}).

If the transformation of (\ref{can:PBC}) to (\ref{desired:PBC:new}) is possible, then the closed-loop system is at least Lyapunov stable with $\widehat H$ as a Lyapunov function since
\begin{align*}
\frac{d\widehat H}{dt} &= \dot q \frac{\partial \widehat H}{\partial q} + \dot p\frac{\partial \widehat H}{\partial p}\\
&= -p^T \widehat M^{-1} G K_{\rm d} G^T \widehat M^{-1}p + C(\widehat M^{-1}p, \widehat M^{-1}p, \widehat M^{-1}p)\\
&= -p^T \widehat M^{-1} G K_{\rm d} G^T \widehat M^{-1}p \leq 0,
\end{align*}
where $C(\widehat M^{-1} p, \widehat M^{-1} p, \widehat M^{-1} p) = 0$ by  (\ref{C:Jacobi:intrinsic})  or statement 1 of Lemma~\ref{first:lemma}. Since the term $\widehat M^{-1}p \righthalfcup \widehat M^{-1} p \righthalfcup  C$ does not change the Hamiltonian $\widehat H$, it is called  a gyroscopic force.

\subsection{Comparison of the New Form of IDA-PBC with the Conventional Form for Mechanical Systems}
\label{subsection:CEI}

In the conventional method of IDA-PBC \cite{AcOrAsMa05,OrGa04}, for a given system of the form (\ref{can:PBC}), the following desired form of controlled Hamiltonian system in place of (\ref{desired:PBC:new}) is used
\begin{align}\label{desired:PBC:old}
\begin{bmatrix}
\dot q \\ \dot p
\end{bmatrix}  = \begin{bmatrix}
 0 & M^{-1}(q)\widehat M(q)\\
 -\widehat M(q)M^{-1}(q) &  J(q,p) -G(q) K_{\rm d}(q) G^T  (q)
 \end{bmatrix} \begin{bmatrix}
\curlyd_q \widehat  H \\ \curlyd_p \widehat H
\end{bmatrix}
\end{align}
with the Hamiltonian (\ref{hat:H}), where $J (q,p)$ is an $n\times n$ skew-symmetric matrix that is linear in $p$. In other words, the $(i, j)$th entry $ J_{ij}(q,p)$ of $J(q,p)$ is written as
\[
 J_{ij}(q,p) =  J_{ij}^k(q)p_k
\]
where $J^k_{ij}$ satisfies
\begin{equation}\label{J:skew}
 J_{ij}^{k} = - J_{ji}^k
\end{equation}
due to the skew symmetry assumption on $J(q,p)$.
Define $B = B_{ijk}(q) e^i \otimes e^j \otimes e^k \in \btensor (\mathbb R^n)$ by
\begin{equation}\label{BJ:relation}
B_{kij} =  J_{ji}^\ell \widehat M_{\ell k},
\end{equation}
which is indeed in $\btensor (\mathbb R^n)$ by (\ref{J:skew}).
Since the matrix $\widehat M$ is invertible, the relationship in (\ref{BJ:relation}) between $(B_{ijk})$ and $(J^k_{ij})$ is bijective.

Let $F^J= F^J_i e^i$ be the force term  $J(q,p) \partial_p\widehat H$ in (\ref{desired:PBC:old}), which can be written component-wise as follows:
\begin{align}
F^J_i &= J^k_{ij} \widehat M^{j\ell}p_k p_\ell \nonumber \\
&= B_{kji}\widehat M^{ks}\widehat M^{j \ell} p_s p_\ell \nonumber \\
&= \frac{1}{2} (B_{kji} + B_{jki})\widehat M^{ks}\widehat M^{j \ell} p_sp_\ell, \label{B:symmetrize}
\end{align}
where the symmetrization  in (\ref{B:symmetrize}) with respect to the indices $k$ and $j$ is valid since $F^J_i$ is a quadratic function of $\widehat M^{-1}p$.  Observe in (\ref{B:symmetrize})  that it is the symmetric part $\frac{1}{2} (B_{kji} + B_{jki})$ of $B_{kji}$ with respect to the first two indices that essentially contributes to the dynamics (\ref{desired:PBC:old}), but  the  skew-symmetric part $\frac{1}{2} (B_{kji} - B_{jki})$ of $B_{kji}$ with respect to the first two indices do not appear in the dynamics at all.  Since $(B_{ijk})$  and $(J^k_{ij})$ are isomorphically related by (\ref{BJ:relation}), it follows that only part of $(J^k_{ij})$ contributes to the dynamics (\ref{desired:PBC:old}).

The natural question that now arises is how to express the essential part of $(J^k_{ij})$ that appears in the dynamics (\ref{desired:PBC:old}), removing all the unnecessary components.
Define $C = C_{ijk}e^i \otimes e^j \otimes e^k \in \ctensor (\mathbb R^n)$ by
\[
C(u,v,w) = \frac{1}{2} (B(u,v,w) + B(v,u,w)) \quad \forall u,v,w \in \mathbb R^n,
\]
or in coordinates
\begin{equation}\label{C:half:B}
C_{ijk} = \frac{1}{2} (B_{ijk} + B_{jik}),
\end{equation}
where $B= B_{ijk}e^i \otimes e^j \otimes e^k \in \btensor (\mathbb R^n)$ is the (0, 3)-tensor field defined in (\ref{BJ:relation}). By Lemma \ref{lemma:BC}, the tensor field $C$ defined above indeed belongs to $\ctensor (\mathbb R^n)$, and such a $C$ expresses exactly the symmetric part $\frac{1}{2} (B_{kji} + B_{jki})$ that contributes to the dynamics (\ref{desired:PBC:old}).
Moreover, from Lemma \ref{lemma:BC}, we can see that the use of the gyroscopic force $\widehat M^{-1} p \righthalfcup \widehat M^{-1} p \righthalfcup  C$ with $C \in \ctensor (\mathbb R^n)$ in (\ref{desired:PBC:new}) contains a smaller number of components than $(J^k_{ij})$ but still makes the same contribution to the dynamics as the term $J(q,p)\widehat M^{-1}p$ does. Concretely speaking, $(J^k_{ij})$ or equivalently $(B_{ijk})$ contains $\frac{n(n-1)(n-2)}{6}$ unnecessary components\footnote{Notice that $\frac{n(n-1)(n-2)}{6} >0$ if and only if the number of degrees of freedom $n$ is greater than 2.} that contribute no terms to the dynamics  (\ref{desired:PBC:old}), showing the superiority of the new form of IDA-PBC in  (\ref{desired:PBC:new}) to the conventional one in (\ref{desired:PBC:old}).
 Hence, we shall exclusively use the new form   (\ref{desired:PBC:new}) in the rest of the paper.

\subsection{Matching Conditions}
\label{subsec:MC}

We study the synthesis problem: what are the matching  conditions for $\widehat H$ and $C$ to satisfy in order to transform (\ref{can:PBC}) to (\ref{desired:PBC:new})?

Let $G^\perp (q)$ denote an $(n-m) \times n$ matrix  whose rows span the left annihilator of the column space of $G(q)$, i.e. $\operatorname{rank} G^\perp = n-m$ and $G^\perp G = 0$.  For the sake of simplicity, we also denote by $G^\perp(q)$ the left annihilator  of $G(q)$.
  Then, by comparing  equations (\ref{can:PBC}) and (\ref{desired:PBC:new}) and collecting terms of equal degrees in $p$, we  obtain   the matching conditions:
\begin{align}
0&= G^\perp ( \curlyd_q V  - \widehat MM^{-1}\curlyd_q\widehat V)\label{potential:matching}\\
0&= G^\perp ( \curlyd_q (p^T M^{-1}p) - \widehat MM^{-1} \curlyd_q (p^T \widehat M^{-1}p)+ 2\widehat M^{-1} p \righthalfcup \widehat M^{-1} p \righthalfcup C  ),\label{kinetic:matching}
\end{align}
where the first set is called the potential matching conditions and the second the kinetic matching conditions. The difference between the two dynamics (\ref{can:PBC}) and (\ref{desired:PBC:new}) is taken care of by the feedback control:
\begin{equation}\label{result:feedback}
u = (G^TG)^{-1} G^T  (\curlyd_q H - \widehat MM^{-1}  \curlyd_q \widehat H  -  GK_{\rm d} G^T \curlyd_p\widehat H + \widehat M^{-1} p \righthalfcup \widehat M^{-1} p \righthalfcup C ).
\end{equation}

The number of PDEs in the kinetic matching conditions would be $\frac{n(n+1)(n-m)}{2}$ if we simply set the coefficients  of  $p_ip_j$'s to zero in (\ref{kinetic:matching}). However, we can reduce the number of PDEs in the kinetic matching conditions by decomposing them into two groups: one without any entries of $C_{ijk}$ and the other with some entries  of $C_{ijk}$. For this purpose, let us introduce a $(2,1)$-tensor field $A = A^{ij}_{k}(q) e_i\otimes e_j \otimes e^k$ and a $(0,3)$-tensor field $\Ttensor = \Ttensor_{ijk}(q) e^i \otimes e^j \otimes e^k$ that are defined as follows:
\begin{equation}\label{intro:A}
A^{ij}_{k} = \frac{1}{2}\widehat M_{k\ell }M^{\ell r}\frac{\partial \widehat M^{ij}}{\partial q^r} - \frac{1}{2}\frac{\partial M^{ij}}{\partial q^k}
\end{equation}
and
\begin{equation}\label{def:T:in:A}
\Ttensor (u,v,w) = A(\widehat Mu, \widehat Mv,w) \quad \forall u,v,w\in \mathbb R^n
\end{equation}
or in coordinates
\begin{align}\label{def:Tijk}
\Ttensor_{ijk} &= A^{rs}_{k}\widehat M_{ri}\widehat M_{sj}\nonumber \\
&= \frac{1}{2}\widehat M_{k\ell }M^{\ell t}\frac{\partial \widehat M^{rs}}{\partial q^t} \widehat M_{ri} \widehat M_{sj}  - \frac{1}{2}\frac{\partial M^{rs}}{\partial q^k} \widehat M_{ri} \widehat M_{sj}\nonumber \\
&= -\frac{1}{2}\widehat M_{k\ell }M^{\ell t} \frac{\partial \widehat M_{ij}}{\partial q^t} - \frac{1}{2}\frac{\partial M^{rs}}{\partial q^k}\widehat M_{ri}\widehat M_{sj}.
\end{align}
Notice that
\[
A^{ij}_{k} = A^{ji}_{k} \quad \textup{and} \quad \Ttensor_{ijk} = \Ttensor_{jik}.
\]
Then, the kinetic matching conditions in (\ref{kinetic:matching}) can be written in the following compact form:
\begin{align}
&C(\widehat M^{-1}p, \widehat M^{-1}p, w) = A(p,p,w) \,\,\, \forall p \in \mathbb R^n, w \in G^\perp\nonumber\\
\Leftrightarrow\, & C(u, u, w) = A(\widehat M u,\widehat M u,w) \,\,\, \forall u \in \mathbb R^n, w \in G^\perp\nonumber \\
\Leftrightarrow\, & C(u, u, w) = \Ttensor(u,u,w)\,\,\, \forall u \in \mathbb R^n, w \in G^\perp\label{holo:cond} \\
\Leftrightarrow\, & C(u, v, w) = \Ttensor(u,v,w) \,\,\, \forall u, v \in \mathbb R^n, w \in G^\perp \label{last:equiv}
\end{align}
where  the equivalence in (\ref{holo:cond}) comes from (\ref{def:T:in:A}), and   the equivalence in (\ref{last:equiv}) is due to  the symmetry  in their first two indices of both $C$ and $\Ttensor$ and the polarization technique, i.e, $4C(u,v,w) = C(u + v, u + v,w) - C(u-v,u-v,w)$. The kinetic matching conditions (\ref{last:equiv}) can be regarded as equations defining the tensor field $C$ on $\mathbb R^n \times \mathbb R^n \times G^\perp$ in terms of $\widehat M$, $M$ and their first-order partial derivatives. By Lemma  \ref{lemma:partial:gyro},  there exists a $(0,3)$-tensor field $C =
 C_{ijk}(q) e^i \otimes e^j \otimes e^k$ satisfying (\ref{C:symmetry:intrinsic}), (\ref{C:Jacobi:intrinsic}) and (\ref{last:equiv}) if and only if   $\Ttensor$   satisfies
\begin{equation}\label{T:Jacobi}
\Ttensor (u,v,w) + \Ttensor (v,w,u) + \Ttensor (w,u,v) = 0 \quad \forall u,v,w\in G^\perp,
\end{equation}
which is a set of PDEs for  $\widehat M$ not containing any components of $C$. The number of PDEs in (\ref{T:Jacobi}) is
\begin{equation}\label{num:PDE:K}
\frac{(n-m+2)(n-m+1)(n-m)}{6}
\end{equation}
which is smaller than $\frac{n(n+1)(n-m)}{2}$ that is  the number of PDEs  in (\ref{kinetic:matching}) before their re-grouping.   Once the PDEs  (\ref{T:Jacobi}) are solved for $\widehat M$,  the tensor field $C$ can be determined purely algebraically as in the proof of Lemma~\ref{lemma:partial:gyro}. 

This approach is an improvement of the conventional IDA-PBC method in \cite{AcOrAsMa05,OrGa04}. In the earlier works,  the decomposition of the kinetic energy into the two groups was never done, so they often  worked with the $\frac{n(n+1)(n-m)}{2}$ PDEs  in (\ref{kinetic:matching}) regarding $J(q,p)$ in (\ref{desired:PBC:old}) (or equivalently  $C$  in the new form  (\ref{desired:PBC:new}))  as a free parameter. For the purpose of comparison, let us consider the case of one degree of underactuation, i.e., $n-m=1$. In the conventional IDA-PBC procedure there are normally more than one PDEs for $\widehat M$ coming from the kinetic matching conditions (\ref{kinetic:matching}); see \cite{AcOrAsMa05} for example.  In our method, however,  the number of PDEs for $\widehat M$ is always one by (\ref{num:PDE:K}).
For example, if $G(q)$ is spanned by $\{e^2, \ldots, e^n\}$, then our kinetic matching condition in (\ref{T:Jacobi}) becomes
\begin{equation}\label{striking:PDE}
\widehat M_{1\ell}M^{\ell r} \frac{\partial \widehat M_{11}}{\partial q^r} +\frac{\partial M^{rs}}{\partial q^1} \widehat M_{r1} \widehat M_{s1} = 0.
\end{equation}
This is only one quasi-linear  PDE for only $ \widehat M_{11}$.  Compare this with the result in \cite{AcOrAsMa05}, where the number of kinetic matching PDEs increases as the number of underactuation degree increases.   This shows the superiority of the new reduced kinetic matching conditions (\ref{T:Jacobi}) that comes from the decomposition of the kinetic matching condition.

From the discussions in Section  \ref{subsection:CEI},  it follows that  (\ref{T:Jacobi})  could be also obtained with the conventional form of IDA-PBC (\ref{desired:PBC:old}), but it has never been done anyway.  Even with this conventional approach, one would have to go through the quadratic gyroscopic force via equations (\ref{BJ:relation}) -- (\ref{C:half:B})  and Lemmas \ref{lemma:BC} and \ref{lemma:partial:gyro}, directly or indirectly. By introducing the quadratic gyroscopic force from the outset, we  avoid  this unnecessary detour.

\subsection{Matching Conditions for Integrable $G(q)$}
\label{subsec:MCIG}

In Section  \ref{subsec:MC} we studied how to reduce the number of PDEs in the kinetic matching condition, by eliminating the gyroscopic term $C$. We now show that the number of partial derivatives of $\widehat M_{ij}$ that appear in the kinetic matching PDEs, can be reduced when the co-distribution generated by the column vectors of  $G(q)$ in the dynamics (\ref{can:PBC}) is integrable, by which we mean the integrability of $G(q)$.

Suppose that $G(q) = \operatorname{span}\{e^{n-m+1}, \ldots, e^n\}$ for each $q$ such that $G^\perp(q) = \operatorname{span}\{e_1, \ldots, e_{n-m} \}$. Let us use the following three groups of indices: $i,j,k,\ldots =1, \ldots, n$; $\alpha, \beta, \gamma, \ldots = 1, \ldots, n-m$; and $a, b, c, \ldots  =n-m+1, \ldots, n$. Then, the potential matching PDEs  (\ref{potential:matching}) are written as
\begin{equation}\label{pot:PDE:coord}
\frac{\partial V}{\partial q^\alpha }  - \widehat M_{\alpha i}M^{ij}\frac{\partial \widehat V}{\partial q^j} = 0.
\end{equation}
The new kinetic matching PDEs
 (\ref{T:Jacobi}) are written as
 \begin{equation}\label{kin:PDE:coord}
 \Ttensor_{\alpha \beta \gamma} + \Ttensor_{ \beta \gamma \alpha} +  \Ttensor_{  \gamma \alpha \beta} = 0
 \end{equation}
 where
 \begin{equation}\label{S:tensorLdef}
 \Ttensor_{\alpha \beta \gamma} = -\frac{1}{2}\widehat M_{\gamma \ell }M^{\ell t} \frac{\partial \widehat M_{\alpha \beta}}{\partial q^t} - \frac{1}{2}\frac{\partial M^{rs}}{\partial q^\gamma}\widehat M_{r\alpha }\widehat M_{s\beta}.
 \end{equation}
Notice that   the partial derivatives of only $\widehat M_{\alpha \beta}$ appear in this set of kinetic matching PDEs (\ref{kin:PDE:coord}) whereas the partial derivatives of $\widehat M_{\alpha a}$ or $\widehat M_{ab}$ do not appear.  It is also remarkable that $\widehat M_{ab}$'s do not appear at all in the matching PDEs in (\ref{pot:PDE:coord}) or (\ref{kin:PDE:coord}).

We now summarize the design procedure for IDA-PBC in the case when $G(q) = \operatorname{span}\{e^{n-m+1}, \ldots, e^n\}$.
\begin{itemize}
\item [P1.] Solve  the  matching PDEs in (\ref{pot:PDE:coord}) and (\ref{kin:PDE:coord}) for $\widehat M (q)$ and $\widehat V (q)$ such that $\widehat M(q)$ is positive definite and $\widehat V(q)$ obtains a non-degenerate minimum at $q=0$. (For local positive definiteness of $\widehat M(q)$, it suffices to have positive definiteness of $\widehat M(0)$)

\item [P2.] Compute $C$ following (\ref{C1}) -- (\ref{C4}), using $\Ttensor_{ijk}$ in (\ref{def:Tijk}).

\item [P3.] Choose an arbitrary $m\times m$ positive definite symmetric matrix $K_{\rm d}$.

\item [P4.] Compute the feedback law  (\ref{result:feedback}).

\end{itemize}
This procedure can be easily adapted to the case of non-integrable $G(q)$.

\subsection{Criteria for Energy Shaping}
We consider the following problem: when does the design procedure in Section  \ref{subsec:MCIG} produce a (locally)  positive definite symmetric matrix $\widehat M$ and a potential function $\widehat V$ with a non-degenerate minimum at $q=0$? In general this is a difficult problem. However, complete answers are available for two cases: 1. when the controlled Hamiltonian is linear and 2. when the degree of underactuation is one, i.e, $n-m=1$. The linear case is  important not only by itself but also for the nonlinear case because the linear result is used in choosing an initial condition for the matching PDEs for the nonlinear case so that $\widehat H$ obtains a non-degenerate minimum at the equilibrium point.
 If the degree of underactuation is more than one, then solving the matching PDEs becomes a challenging job, which would normally involve the operations of prolongation and projection on the system of matching PDEs; refer to  \cite{Po94} for definitions of prolongation and projection. Luckily, when the degree of underactuation is one, no knowledge of  formal theory of systems of PDEs  is necessary, but a simple application of  Frobenius' integrability theorem  suffices.
 
 \subsubsection{Case of Linear Controlled Hamiltonian Systems}

Let us first take the linear case. Recall that the controlled Hamiltonian system (\ref{can:PBC})  is called linear if $M$ is a constant matrix, $V$ is a quadratic function of $q$, and $G$ is a constant matrix. For  IDA-PBC of linear controlled Hamiltonian systems, the gyroscopic term $C$ in (\ref{desired:PBC:new}) is not necessary since it would produce  nonlinear terms in the dynamics, so we set $C=0$. More concretely, let 
\[
H(q,p) = \frac{1}{2}p^TM^{-1}p + \frac{1}{2}q^TSq, \quad \widehat H(q,p) = \frac{1}{2}p^T \widehat M^{-1}p + \frac{1}{2}q^T\widehat Sq
\]
where 
\[
M = M^T \succ 0, \quad S=S^T
\]
and
\[
 \widehat M = \widehat M^T \succ 0, \quad \widehat S = \widehat S^T \succ 0. 
\]
Here, $\widehat M$ and $\widehat S$ are to be found. Then, the original  controlled Hamiltonian system (\ref{can:PBC})  is given by
\begin{align}
\dot q &= M^{-1}p,  \label{Loriq}\\
\dot p &= -S q + Gu \label{Lorip}
\end{align}
and the desired form of controlled Hamiltonian system (\ref{desired:PBC:new}) with $C=0$ is given by
\begin{align}
\dot q &= M^{-1}p, \label{dLLq}\\
\dot p &= -\widehat MM^{-1}\widehat S q -GK_{\rm d} G^T \widehat M^{-1}p, \label{dLLp}
\end{align}
where $K_{\rm d} = K^T_{\rm d} \succ 0$. 
Let us choose a feedback control $u$ of the form
\begin{equation}\label{u:lin:cont}
u = -K_0 q - K_{\rm d}G^T \widehat M^{-1} p,
\end{equation}
which transforms the system in  (\ref{Loriq}) and (\ref{Lorip}) to the following closed-loop system:
\begin{align}
\dot q &= M^{-1}p,  \label{Foriq}\\
\dot p &= -(S+GK_0) q -GK_{\rm d} G^T \widehat M^{-1}p. \label{Forip}
\end{align}
For the closed-loop system in (\ref{Foriq}) and (\ref{Forip}) to be equal to the desired system in (\ref{dLLq}) and (\ref{dLLp}), we need to have
\[
\widehat M M^{-1}\widehat S = S+GK_0,
\]
or
\begin{equation}\label{linear:mat}
M^{-1}\widehat M M^{-1}\widehat S = M^{-1}S+M^{-1}GK_0.
\end{equation}
Hence, our problem comes to down finding a matrix $K_0$ such that there exist two symmetric positive definite matrices $\widehat M$ and $\widehat S$ that satisfy (\ref{linear:mat}).

\begin{lemma}\label{lin:Ham:lem}
There is a  matrix $K_0$ such that there exist two symmetric positive definite matrices $\widehat M$ and $\widehat S$ that satisfy (\ref{linear:mat}), if and only if the system in (\ref{Loriq}) and (\ref{Lorip}) is controllable or its uncontrollable part is  oscillatory.\footnote{A linear dynamics $\dot x = Ax$ is called oscillatory if $A$ is diagonalizable and each eigenvalue of $A$ is a non-zero purely imaginary number.}
\begin{proof}
The Lagrangian equivalent of this theorem  is  proved in \cite{Ch10:SIAM} (see also \cite{Ch12:IJRNC} for a direct proof) and it can be easily adapted to the Hamiltonian case. Hence, instead of repeating the proof in \cite{Ch10:SIAM,Ch12:IJRNC} we here give a brief algorithm how to choose such a matrix $K_0$.  Suppose that the system in (\ref{Loriq}) and (\ref{Lorip}) is controllable. Then, one can assign arbitrary eigenvalues to the matrix  $M^{-1}S+M^{-1}GK_0$ by choosing an appropriate matrix $K_0$ provided that complex eigenvalues appear in conjugate pair. Hence, choose a matrix $K_0$ such that the matrix  $M^{-1}S+M^{-1}GK_0$ is diagonalizable and has positive (real) eigenvalues only: for example, we can make    $M^{-1}S+M^{-1}GK_0$ have distinct positive eigenvalues. Let $\lambda_1, \ldots, \lambda_n$ be the eigenvalues of $M^{-1}S+M^{-1}GK_0$. Then, there is an $n\times n$ matrix $J$ such that 
\[
M^{-1}S+M^{-1}GK_0  = J\textup{diag}\{ \lambda_1 \oplus \lambda_2 \oplus \cdots \oplus \lambda_n\}J^{-1}
\]
 with all $\lambda_i >0$.  Choose $n$ arbitrary positive numbers $\mu_1, \ldots, \mu_n$, and let 
 \[
 \widehat{M} = MJ \textup{diag}\{ \mu_1 \oplus \cdots \oplus \mu_n\}J^TM
 \]
  and 
\[
\widehat{S} = (J^{-1})^T \textup{diag}  \left \{ \frac{\lambda_1}{\mu_1} \oplus \cdots \oplus \frac{\lambda_n}{\mu_n} \right \}J^{-1}.
\]
Both $\widehat M$ and $\widehat S$ are symmetric and positive definite, and they satisfy  (\ref{linear:mat}).  The case in which the system in (\ref{Loriq}) and (\ref{Lorip}) is not controllable but its uncontrollable part is oscillatory can be handled similarly; see  \cite{Ch10:SIAM,Ch12:IJRNC}  for more detail. 
\end{proof}
\end{lemma}

Now, suppose that we have found such a matrix $K_0$. Then,   the closed-loop system  in (\ref{dLLq}) and (\ref{dLLp}), or equivalently (\ref{Foriq}) and (\ref{Forip}), is  Lyapunov stable. It turns out that the closed-loop system is exponentially stable if and only if the original system in (\ref{Loriq}) and (\ref{Lorip}) is controllable; see \cite{Ch10:SIAM,Ch12:IJRNC} for a proof of this statement.  The following theorem summarizes the results discussed so far. 

\begin{theorem}\label{thm:linear:Ham}
The controlled Hamiltonian system (\ref{Loriq}) and (\ref{Lorip})  can be transformed to the desired form (\ref{dLLq}) and (\ref{dLLp}) with $\widehat M = \widehat M^T \succ 0$ and $\widehat S=\widehat S^T \succ 0$ if and only if the system (\ref{Loriq}) and (\ref{Lorip})  is controllable or its uncontrollable part is oscillatory. Moreover, the closed-loop system  (\ref{dLLq}) and (\ref{dLLp}) is exponentially stable if and only if the original system  (\ref{Loriq}) and (\ref{Lorip}) is controllable.
\end{theorem}
Notice that the conditions in Theorem \ref{thm:linear:Ham} are on the given system (\ref{Loriq}) and (\ref{Lorip}) and purely algebraic, so they are easy to verify. 
Let us revisit the example in (\ref{int:old:eq}) -- (\ref{intro:H1}). From (\ref{intro:H0})
\[
M = \begin{bmatrix}
1 & 0 \\
0 & 1
\end{bmatrix}, \quad 
S = \begin{bmatrix}
-1 & 1 \\ 1 & 0
\end{bmatrix}.
\]
From (\ref{int:old:eq})
\[
G =\begin{bmatrix}
0 \\ 1
\end{bmatrix}.
\]
The system (\ref{int:old:eq}) is controllable, so energy shaping is possible on the system and the resultant closed-loop system is exponentially stable by Theorem \ref{thm:linear:Ham}. 
In (\ref{intro:u:2}) we have chosen
\[
K_0 = \begin{bmatrix}-7 &  4 \end{bmatrix}
\]
such that 
\[
\lambda_1 = 1, \quad \lambda_2 = 2.
\]
The matrix $J$ in the proof of Lemma \ref{lin:Ham:lem} is given by
\[
J = \begin{bmatrix}
1 & 1\\ 2 & 3
\end{bmatrix}.
\]
We have chosen 
\[
\mu_1 = \mu_2 = 1
\]
such that
\[
\widehat M = \begin{bmatrix}
2 & 5\\
5 & 13
\end{bmatrix}, \quad 
\widehat S = \begin{bmatrix}
17 & -7 \\ -7 & 3
\end{bmatrix}
\]
both of which are symmetric and positive definite.   We set  $K_{\rm d} =1$  so that (\ref{u:lin:cont}) yields (\ref{intro:u:2}). One can easily check that the closed-loop system is indeed exponentially stable.

 \subsubsection{Case of Nonlinear Controlled Hamiltonian Systems with Underactuation Degree One}
We now consider the case of nonlinear controlled Hamiltonian systems with  underactuation degree one, i.e. $n-m = 1$.
\begin{theorem}\label{thm:nonlinear:PBC:ud1}
Consider the controlled Hamiltonian system   (\ref{can:PBC}) with one degree of underactuation. Let $\Sigma^\ell$ denote its linearization at the equilibrium  $(q,p) = (0,0)$. The system (\ref{can:PBC})  can be transformed to the desired form (\ref{desired:PBC:new}) with $\widehat H$ having a non-degenerate local minimum at $(q,p) = (0,0)$, $K_{\rm d} = K_{\rm d}^T$ being positive definite, and $C \in \ctensor (\mathbb R^n)$, if and only if $\Sigma^\ell$ is controllable or the uncontrollable part of $\Sigma^\ell$ is oscillatory. Moreover, the closed-loop system (\ref{desired:PBC:new}) is (locally) exponentially stable if and only if $\Sigma^\ell$ is controllable.
\begin{proof}
Without loss of generality,  we may assume $G(q) = \operatorname{span}\{ e^2, \ldots, e^n\}$. Then, the potential and kinetic matching conditions in (\ref{pot:PDE:coord}) and (\ref{kin:PDE:coord}) can be written as
\begin{align}
 \widehat M_{1i}M^{ij}\frac{\partial \widehat V}{\partial q^j} -\frac{\partial V}{\partial q^1} &= 0 \label{ud1:pot}\\
\widehat M_{1i}M^{ij} \frac{\partial \widehat M_{11}}{\partial q^j} +\frac{\partial M^{ij}}{\partial q^1} \widehat M_{i1} \widehat M_{j1} &= 0. \label{ud1:kin}
\end{align}
Write the Hamiltonian $H^\ell$ of the linearization $\Sigma^\ell$ of (\ref{can:PBC}) at the origin as
\[
H^\ell = \frac{1}{2}p^T M^{-1}(0) p + \frac{1}{2}q^T D^2V(0) q
\]
where $M^{-1}(0)$ is the evaluation of $M^{-1}(q)$ at $q=0$, and $D^2V$ is the second-order derivative matrix of $V$.
Suppose that $\Sigma^\ell$ is controllable or its uncontrollable part is oscillatory. Then, by Theorem~\ref{thm:linear:Ham} there are two constant positive definite symmetric matrices $\overline M $ and $\overline S$ such that
\begin{equation}\label{jinkook}
\overline M_{1i}M^{ij}(0) \overline S_{jk}  - \frac{\partial^2 V}{\partial q^1 \partial q^k}(0) = 0
\end{equation}
which is  the potential matching condition for $\Sigma^\ell$.
The kinetic matching condition trivially holds for $\Sigma^\ell$ since both $M$  and $\overline M$ are  constant matrices. Notice that the set of PDEs in (\ref{ud1:pot}) and (\ref{ud1:kin}) for $\widehat V$ and $\widehat M_{11}$ is integrable, which can be easily checked by applying Frobenius'  integrability theorem; refer to  \cite{Le02} for Frobenius' integrability theorem.  We now impose the following initial conditions on $\widehat M$ and $\widehat V$ at $q=0$:
\begin{equation}\label{init:cond:Hhat}
\widehat M(0) = \overline M, \quad {\mathbf d}\widehat V(0) = 0, \quad D^2\widehat V(0) = \overline S.
\end{equation}
 The initial conditions  at $q=0$ are compatible with the PDE (\ref{ud1:pot}) since the differentiation of (\ref{ud1:pot}) and its evaluation  at $q=0$ yield
 \begin{equation}\label{init:cond:2}
 \widehat M_{1i}(0) M^{ij}(0) \frac{\partial^2 \widehat V}{\partial q^j \partial q^k}(0) -  \frac{\partial^2 V}{\partial q^1 \partial q^k}(0) =0
 \end{equation}
 where ${\mathbf d}\widehat V(0) = 0$ is used. Equation  (\ref{init:cond:2}) has the same structure as  (\ref{jinkook}). Hence, by Frobenius' integrability theorem there exist $\widehat M = \widehat M^T$ and $\widehat V$ that satisfy (\ref{ud1:pot}), (\ref{ud1:kin}) and (\ref{init:cond:Hhat}). Due to the initial conditions (\ref{init:cond:Hhat}), the desired Hamiltonian $\widehat H = \frac{1}{2}p^T\widehat M^{-1} p + \widehat V(q)$ obtains a non-degenerate local minimum value at $(q,p) = (0,0)$. Choose $C$, $K_{\rm d}$ and $u$ following the procedure in Section  \ref{subsec:MCIG}.  We have successfully transformed the given system (\ref{can:PBC}) to the desired form (\ref{desired:PBC:new}).

Suppose now that the given nonlinear controlled Hamiltonian system (\ref{can:PBC}) can be transformed to the desired form (\ref{desired:PBC:new}) with $\widehat H = \frac{1}{2}p^T \widehat M^{-1}p + \widehat V(q)$, where $\widehat H$ has a non-degenerate minimum value at $(q,p) = (0,0)$. Through linearization, we can see that the linearization $\Sigma^\ell$ of (\ref{can:PBC}) gets transformed to the following linear dynamics:
\begin{align}\label{lin:desired:PBC:new}
\begin{bmatrix}
\dot q \\ \dot p
\end{bmatrix}  = \begin{bmatrix}
 0 & M^{-1}(0)\widehat M(0)\\
 -\widehat M(0)M^{-1}(0) & -G(0) K_{\rm d}(0) G^T  (0)
 \end{bmatrix}\begin{bmatrix}
\curlyd_q \widehat  H^\ell \\ \curlyd_p \widehat H^\ell
\end{bmatrix} ,
\end{align}
where  $\widehat H^\ell  = \frac{1}{2}p^T\widehat M^{-1}(0)p + \frac{1}{2}q^T D^2\widehat V(0)q$ with both $\widehat M(0) $ and $D^2\widehat V(0) $ positive definite and symmetric.
By Theorem~\ref{thm:linear:Ham} the linear  controlled Hamiltonian system  $\Sigma^\ell$ is controllable or its uncontrollable dynamics are oscillatory. This completes the proof of the first statement that the system (\ref{can:PBC})  can be transformed to the desired form (\ref{desired:PBC:new}) if and only if $\Sigma^\ell$ is controllable or its uncontrollable part of $\Sigma^\ell$ is oscillatory.

We now prove the statement on exponential stability.  By the Lyapunov linearization theorem, the closed-loop system (\ref{desired:PBC:new}) is exponentially stable if and only if its linearized dynamics, say the one in (\ref{lin:desired:PBC:new}),  is exponentially stable. From Theorem~\ref{thm:linear:Ham}, we know that (\ref{lin:desired:PBC:new}) is exponentially stable if and only if the linear system $\Sigma^\ell$, which is the linearization of the original nonlinear controlled Hamiltonian system (\ref{can:PBC}), is controllable. This completes the proof.

\end{proof}
\end{theorem}

Theorem~\ref{thm:nonlinear:PBC:ud1} is important because it gives necessary and sufficient conditions for Lyapunov/exponential stabilizability by IDA-PBC for the case of underactuation degree one. Moreover,
 the conditions therein are on the given system (\ref{can:PBC}) that can be easily verified in advance before any attempt to find a desired system (\ref{desired:PBC:new}).  In the past, they had only sufficient conditions for stabilizability by IDA-PBC, so they were inconclusive on stabilizability by IDA-PBC when the sufficient conditions did not hold.

 Theorem~\ref{thm:nonlinear:PBC:ud1} applies to a wide range of systems including the inverted pendulum on a car, the Pendubot, the Furuta pendulum, the ball and beam system, and the planar vertical take off and landing aircraft,  the linearization of each of which is controllable; refer to \cite{NgChSo13,LiYu13} for the dynamics of these systems. We  remark that the  feedback control (\ref{result:feedback})  is nonlinear.

\subsection{Inverted Pendulum on a Cart}

\begin{figure}[t]
\begin{center}
\includegraphics[scale = .5]{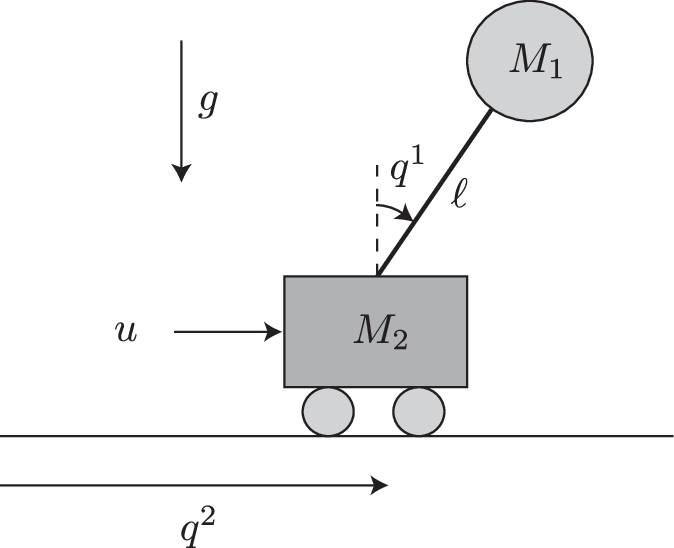}
\caption{\label{ipc2.figure} An inverted pendulum on a cart.}
\end{center}
\end{figure}

We  consider the inverted pendulum on a cart in Figure~\ref{ipc2.figure}. This system has been stabilized by old  methods of energy shaping in \cite{AcOrAsMa05,BlChLeMa01}, and we here stabilize it with the new IDA-PBC framework.
The mass matrix $M$ and the potential function $V$ are given by
\[
M = \begin{bmatrix}
M_1\ell^2 & M_1\ell \cos q^1 \\
M_1 \ell \cos q^1 & M_1 + M_2
\end{bmatrix}
\]
and
\[
V(q) = M_1g\ell (1-\cos q^1).
\]
Since a control $u$  is given in the $q^2$ direction, we have
\[
G(q) = \begin{bmatrix}
0 \\ 1
\end{bmatrix}.
\]
Our main purpose is to illustrate the design procedure in Section  \ref{subsec:MCIG}, so we choose the following values of parameters:
\[
M_1 = M_2 = \ell = 1, \quad  g = 10.
\]
One can easily check that the linearization of this system at the origin is controllable. Hence, by Theorem~\ref{thm:nonlinear:PBC:ud1} we can exponentially stabilize this system with the  method of IDA-PBC.

Let us now construct a control law following the procedure in Section  \ref{subsec:MCIG}.
The  potential and kinetic matching conditions in (\ref{pot:PDE:coord}) and (\ref{kin:PDE:coord}) are given, after simplification, as
\[
(2\widehat M_{11} - \widehat M_{12} \cos q^1) \frac{\partial \widehat V}{\partial q^1} + (-\widehat M_{11} \cos q^1 + \widehat M_{12} ) \frac{\partial \widehat V}{\partial q^1} = 10\sin q^1( 2 - \cos^2 q^1),
\]
and
\begin{align*}
(2\widehat M_{11} - \widehat M_{12} \cos q^1)& \frac{\partial \widehat M_{11}}{\partial q^1} + (-\widehat M_{11} \cos q^1 + \widehat M_{12} ) \frac{\partial \widehat M_{11}}{\partial q^1}
 \\&= \frac{2\sin q^1(2\widehat M_{11} - \widehat M_{12} \cos q^1) (\widehat M_{12} - \widehat M_{11}  \cos q^1)}{2-\cos^2 q^1}.
\end{align*}
A solution is found to be
\begin{align*}
&\widehat M_{11} = 2\cos^2 q^1-\epsilon, \qquad \widehat M_{12} = (4-\epsilon) \cos q^1,\\
&\widehat V = -\frac{10}{\epsilon}\cos q^1 + \left (q^2 + \frac{2\sin q^1}{\epsilon} \right )^2
\end{align*}
where $0<\epsilon<2$ is a constant. To make $\widehat M$ positive definite at least locally around $q=0$, we choose an $\widehat M_{22}$ such that $\widehat M_{22} > (\widehat M_{12})^2/ \widehat M_{11}$ at $q=0$. For example,
\[
\widehat M_{22} =K+ \frac{(4-\epsilon)^2 \cos^2 q^1}{2\cos^2 q^1 - \epsilon},
\]
where $K>0$ is a constant parameter. It is easy to see that $(q,p) = (0,0)$ is the minimum point of the Hamiltonian $\widehat H(q,p) = \frac{1}{2}p^T\widehat M^{-1}p + \widehat V(q)$ over the set $I_\epsilon \times \mathbb R \times \mathbb R^2$, where
\[
I_\epsilon := \cos^{-1}\left (\left (\sqrt{\frac{\epsilon}{2}}, 1\right ] \right ) \subset \left (-\frac{\pi}{2}, \frac{\pi}{2} \right ).
\]
Notice that $I_\epsilon \rightarrow (-\frac{\pi}{2}, \frac{\pi}{2})$ as $\epsilon \rightarrow 0$.
Compute the gyroscopic tensor $C$ using (\ref{C1}) -- (\ref{C4}) and  (\ref{def:Tijk}) as follows:
\begin{align*}
C_{111} &= C_{222} = 0,\\
C_{121} &= C_{211} = S_{121}\\
C_{221} &=     S_{221}  \\
C_{112} &= -2 S_{121} \\
C_{122} & = C_{212} = -\frac{1}{2}S_{221},
\end{align*}
where
\begin{align*}
S_{121} &= -\frac{1}{2}\left (\widehat M_{1r}M^{rs}\frac{\partial \widehat M_{12}}{\partial q^s}  + \frac{\partial M^{rs}}{\partial q^1}\widehat{M}_{r1}\widehat{M}_{s2} \right ),\\
S_{221} &=  -\frac{1}{2}\left (\widehat M_{1r}M^{rs}\frac{\partial \widehat M_{22}}{\partial q^s}  + \frac{\partial M^{rs}}{\partial q^1}\widehat{M}_{r2}\widehat{M}_{s2} \right ).
\end{align*}
 Choose a $K_{\rm d} =K_{\rm d}^T \succ 0$ and compute the feedback control law (\ref{result:feedback}) which simplifies to
 \[
 u = -\widehat M_{2i}M^{ij} \left ( \frac{1}{2}\frac{\partial \widehat M^{k\ell}}{\partial q^j}p_kp_\ell + \frac{\partial \widehat V}{\partial q^j}\right ) - K_{\rm d}\widehat M^{2i}p_i + C_{ij2}\widehat M^{ik}\widehat M^{j\ell}p_k p_\ell.
 \]
  The closed-loop system is exponentially stable by Theorem~\ref{thm:nonlinear:PBC:ud1}.

\subsection{The Duality of the Method of IDA-PBC to the Method of Controlled Lagrangians}

We show that   the method of IDA-PBC and the method of controlled Lagrangians (CLs) are dual to each other in some sense. Both methods are  for energy shaping and the only difference between the two is that the CL method employs the position and velocity variables $(q,\dot q)$ whereas the IDA-PBC method uses the position and momentum variables $(q,p)$; refer to \cite{Ch10:SIAM,Ch10:TAC,NgChLa13,Ch12:IJRNC} for theory of the CL method and its generalization to non-Lagrangian second-order systems. We here show that the set of matching conditions that arise in the CL method is the same as that in the IDA-PBC method. 

In the CL method the matching conditions involved are equations (3.23) and (3.25), or equations (3.28) and (3.30),  in  \cite{Ch10:SIAM}, which are re-written here as follows:
\begin{align}
&\frac{\partial V}{\partial q^\alpha} - \widehat T_{\alpha i}M^{ij} \frac{\partial \widehat V}{\partial q^j}=0,\label{CL:mat:1}\\
&\widehat J_{\alpha\beta\gamma} + \widehat J_{\beta\gamma\alpha}  + \widehat J_{\gamma\alpha\beta}  = 0, \label{CL:mat:2}
\end{align}
where
\begin{align*}
\widehat{J}_{\alpha\beta\gamma} &= \frac{1}{2} \widehat{T}_{\gamma s}M^{sk} \left (\frac{\partial \widehat{T}_{\alpha \beta}}{\partial q^k} - \widehat{T}_{\alpha i}\Gamma^i_{\beta k} -   \widehat{T}_{\beta i} \Gamma^i_{\alpha k}\right )
\end{align*}
and
\[
\widehat T = M \widehat M^{-1}M, \quad \widehat T_{ij} = M_{ir}\widehat M^{rs} M_{sj}.
\]
Here, $\Gamma^k_{ij}$'s are the Christoffel symbols of the second kind of the mass matrix $M = (M_{ij})$:
\[
\Gamma^k_{ij} = \frac{1}{2} M^{k\ell}\left (  \frac{\partial M_{\ell j}}{\partial q^i} + \frac{\partial M_{i\ell}}{\partial q^j} - \frac{\partial M_{ij}}{\partial q^\ell}\right ).
\]

\begin{theorem}\label{equiv:thm:CL:IDA}
The matching PDEs in (\ref{CL:mat:1}) and (\ref{CL:mat:2}) for $\widehat T_{ij}$ and $\widehat V$ are the same as the PDEs in (\ref{pot:PDE:coord}) and (\ref{kin:PDE:coord}) for $\widehat M_{ij}$ and $\widehat V$. In other words, if one simply replaces $\widehat T_{ij}$ with $\widehat M_{ij}$ in (\ref{CL:mat:1}) and (\ref{CL:mat:2}), then (\ref{pot:PDE:coord}) and (\ref{kin:PDE:coord}) are obtained.
\begin{proof}
It is easy to see that  the PDE in (\ref{CL:mat:1})  becomes the PDE in  (\ref{pot:PDE:coord}) if $\widehat T_{ij}$'s are replaced with $\widehat M_{ij}$'s. 
For notational convenience, for any tensor $R_{\alpha \beta \gamma}$  define
\[
\sum_{\alpha, \beta, \gamma,  \textup{cyclic}} R_{\alpha \beta \gamma} := R_{\alpha\beta\gamma} + R_{\beta\gamma\alpha}  +R_{\gamma\alpha\beta} .
\]
Notice
\begin{align*}
\sum_{\alpha, \beta, \gamma,  \textup{cyclic}} & \widehat{T}_{\gamma s}M^{sk} \left ( \widehat{T}_{\alpha i}\Gamma^i_{\beta k} +   \widehat{T}_{\beta i} \Gamma^i_{\alpha k}\right ) \\
&=\sum_{\alpha, \beta, \gamma,  \textup{cyclic}}  \frac{1}{2}\widehat T_{\gamma s} \widehat T_{\alpha i} M^{sk}M^{i\ell} \left ( \frac{\partial M_{\ell k}}{\partial q^\beta} + \frac{\partial M_{\beta \ell}}{\partial q^k} - \frac{\partial M_{\beta k}}{\partial q^\ell}\right ) \\
&\quad+ \sum_{\alpha, \beta, \gamma,  \textup{cyclic}}  \frac{1}{2}\widehat T_{\gamma s} \widehat T_{\beta i} M^{sk}M^{i\ell} \left ( \frac{\partial M_{\ell k}}{\partial q^\alpha} + \frac{\partial M_{\alpha \ell}}{\partial q^k} - \frac{\partial M_{\alpha k}}{\partial q^\ell}\right ) \\
&= \sum_{\alpha, \beta, \gamma,  \textup{cyclic}}  \frac{1}{2}\widehat T_{\gamma s} \widehat T_{\alpha i} M^{sk}M^{i\ell}  \frac{\partial M_{\ell k}}{\partial q^\beta} + \sum_{\alpha, \beta, \gamma,  \textup{cyclic}}  \frac{1}{2}\widehat T_{\gamma s} \widehat T_{\beta i} M^{sk}M^{i\ell}  \frac{\partial M_{\ell k}}{\partial q^\alpha} \\
&\quad+\sum_{\alpha, \beta, \gamma,  \textup{cyclic}}  \frac{1}{2}\widehat T_{\gamma s} \widehat T_{\alpha i} M^{sk}M^{i\ell} \left ( \frac{\partial M_{\beta \ell}}{\partial q^k} - \frac{\partial M_{\beta k}}{\partial q^\ell}\right )\\
&\quad+\sum_{\alpha, \beta, \gamma,  \textup{cyclic}}  \frac{1}{2}\widehat T_{\gamma s} \widehat T_{\beta i} M^{sk}M^{i\ell} \left (  \frac{\partial M_{\alpha \ell}}{\partial q^k} - \frac{\partial M_{\alpha k}}{\partial q^\ell}\right ) \\
&= \sum_{\alpha, \beta, \gamma,  \textup{cyclic}}  \frac{1}{2}\widehat T_{\gamma s} \widehat T_{\alpha i} M^{sk}M^{i\ell}  \frac{\partial M_{\ell k}}{\partial q^\beta} + \sum_{\alpha, \beta, \gamma,  \textup{cyclic}}  \frac{1}{2}\widehat T_{\alpha s} \widehat T_{\gamma i} M^{sk}M^{i\ell}  \frac{\partial M_{\ell k}}{\partial q^\beta} \\
&\quad+\sum_{\alpha, \beta, \gamma,  \textup{cyclic}}  \frac{1}{2}\widehat T_{\gamma s} \widehat T_{\alpha i} M^{sk}M^{i\ell} \left ( \frac{\partial M_{\beta \ell}}{\partial q^k} - \frac{\partial M_{\beta k}}{\partial q^\ell}\right )\\
&\quad+\sum_{\alpha, \beta, \gamma,  \textup{cyclic}}  \frac{1}{2}\widehat T_{\alpha s} \widehat T_{\gamma i} M^{sk}M^{i\ell} \left (  \frac{\partial M_{\beta \ell}}{\partial q^k} - \frac{\partial M_{\beta k}}{\partial q^\ell}\right )\\
&= \sum_{\alpha, \beta, \gamma,  \textup{cyclic}} \widehat T_{\gamma s} \widehat T_{\alpha i} M^{sk}M^{i\ell}  \frac{\partial M_{\ell k}}{\partial q^\beta} \\
&=\sum_{\alpha, \beta, \gamma,  \textup{cyclic}}  -\widehat T_{\gamma s} \widehat T_{\alpha i} \frac{\partial M^{is}}{\partial q^\beta}.
\end{align*}
Hence,
\begin{align*}
\sum_{\alpha, \beta, \gamma,  \textup{cyclic}} \widehat J_{\alpha \beta \gamma}  &= \sum_{\alpha, \beta, \gamma,  \textup{cyclic}} \frac{1}{2} \widehat{T}_{\gamma s}M^{sk} \left (\frac{\partial \widehat{T}_{\alpha \beta}}{\partial q^k} - \widehat{T}_{\alpha i}\Gamma^i_{\beta k} -   \widehat{T}_{\beta i} \Gamma^i_{\alpha k}\right ) \\
&= \sum_{\alpha, \beta, \gamma,  \textup{cyclic}} \frac{1}{2} \widehat{T}_{\gamma s}M^{sk} \frac{\partial \widehat{T}_{\alpha \beta}}{\partial q^k} +
\sum_{\alpha, \beta, \gamma,  \textup{cyclic}} \frac{1}{2} \widehat T_{\gamma s} \widehat T_{\alpha i} \frac{\partial M^{is}}{\partial q^\beta} \\
&= \sum_{\alpha, \beta, \gamma,  \textup{cyclic}} \frac{1}{2} \widehat{T}_{\gamma s}M^{sk} \frac{\partial \widehat{T}_{\alpha \beta}}{\partial q^k} +
\sum_{\alpha, \beta, \gamma,  \textup{cyclic}} \frac{1}{2} \widehat T_{\alpha s} \widehat T_{\beta i} \frac{\partial M^{is}}{\partial q^\gamma} \\
&= \sum_{\alpha, \beta, \gamma,  \textup{cyclic}} \left ( \frac{1}{2} \widehat{T}_{\gamma \ell }M^{\ell t} \frac{\partial \widehat{T}_{\alpha \beta}}{\partial q^t}   +\frac{1}{2} \frac{\partial M^{rs}}{\partial q^\gamma} \widehat T_{r \alpha }\widehat T_{s\beta} \right ).  
\end{align*}
It is now easy to see that the PDE in (\ref{CL:mat:2})  becomes the PDE in  (\ref{kin:PDE:coord}) if $\widehat T_{ij}$'s are replaced with $\widehat M_{ij}$'s. 

\end{proof}
\end{theorem}

From Theorem \ref{equiv:thm:CL:IDA}, we see that the CL method is ``dual'' to the IDA-PBC method with the following correspondence between new mass matrices: 
\[
\begin{array}{ccc}
\textup{CL method}&\longleftrightarrow  &\textup{IDA-PBC method}\\
\widehat M & \longleftrightarrow & M \widehat M^{-1}M \\
M\widehat M^{-1}M &     \longleftrightarrow    & \widehat M.
\end{array}
\]

\section{Conclusions and Future Work}

We have improved the method of IDA-PBC in several ways. First, we have showed that there is redundancy in the skew-symmetric interconnection matrix in the conventional form of IDA-PBC, and  then have replaced the skew-symmetric interconnection matrix  term with a gyroscopic force to remove the redundancy. We  have derived the matching conditions and  decomposed them into two parts. As a result of the decomposition, we have come up with a smaller number of kinetic matching PDEs than those in the literature. Moreover,   a smaller number of entries of the desired mass matrix appear in the kinetic matching PDEs when  controls are given in coordinate directions.  Easily verifiable necessary and sufficient conditions have been given for Lyapunov/exponential stabilizability by IDA-PBC for all linear controlled Hamiltonian (or simple mechanical) systems with  arbitrary degree of underactuation and for all nonlinear controlled Hamiltonian (or simple mechanical) systems with one degree of underactuation.  A step-by-step IDA-PBC synthesis procedure was provided and   illustrated with examples. The duality of the IDA-PBC method to the controlled Lagrangian method was discussed.    It will be interesting to extend the energy shaping method to locomotion, tracking  and hybrid mechanical systems \cite{GrBrSp10,HoLeSp07,SpHo07} and to combine it with the recent resents in   \cite{ReRi13} on equivariants of mechanical control systems.

\begin{acknowledgments}
The author would like to thank the late professor Jerrold E. Marsden for his ever-lasting lessons. 
\end{acknowledgments}

\end{document}